\documentclass[a4paper,american,reqno]{amsart}

\newif\ifPreprint \Preprintfalse
\newif\ifSubmission \Submissiontrue

\usepackage{babel}
\usepackage[utf8]{inputenc}
\usepackage[T1]{fontenc}
\usepackage{lmodern}
\usepackage{xspace}
\usepackage{mathtools}
\usepackage{amsmath}
\usepackage{enumitem}
\usepackage{accents}
\usepackage{todonotes}
\usepackage{csquotes}
\usepackage{microtype}
\usepackage{bbm}
\usepackage[style = numeric-comp,
            firstinits = true,
            maxbibnames = 100,
            maxcitenames = 2,
            isbn = false,
            backend = biber]{biblatex}
\usepackage[colorlinks,
            citecolor=blue,
            urlcolor=blue,
            linkcolor=blue]{hyperref}
\usepackage{cleveref}

\makeatletter
\patchcmd{\@settitle}{\uppercasenonmath\@title}{\scshape\large}{}{}
\patchcmd{\@setauthors}{\MakeUppercase}{\scshape\normalsize}{}{}
\makeatother

\vfuzz \hfuzz
\theoremstyle{plain}

\newtheorem{coroll}{Corollary}
\newtheorem{theorem}{Theorem}

\theoremstyle{definition}

\newtheorem{example}{Example}
\newtheorem{definition}{Definition}

\newtheorem{prop}{Proposition}

\theoremstyle{remark}

\newcommand{\defsep}{\colon}
\newcommand{\Unc}{\mathcal{U}}
\newcommand{\Mprod}[1]{\langle #1\rangle}

\newcommand{\X}{\mathcal{X}}
\newcommand{\R}{\mathbb{R}}

\newcommand{\N}{\mathbb{N}}

\newcommand{\fa}{\text{ for all }}
\newcommand{\Zc}{\mathcal{Z}}

\newcommand{\Scal}{\mathcal{S}}

\newcommand{\XRO}{\mathcal{X}^{\mathrm{RO}}}
\newcommand{\XPRO}{\mathcal{X}^{\mathrm{PRO}}}

\newcommand{\rev}[1]{\textcolor{black}{#1}}

\bibliography{bibliography}

\usepackage{amsmath}
\usepackage{amssymb}

\begin{document}

\title[Pareto Robust Optimization on Euclidean Vector Spaces]%
{Pareto Robust Optimization on Euclidean Vector Spaces}
\author[D. Adelhütte, C. Biefel, M. Kuchlbauer, J. Rolfes]%
{Dennis Adelhütte, Christian Biefel, Martina Kuchlbauer, Jan Rolfes}

\address[D. Adelhütte, C. Biefel, M. Kuchlbauer, J. Rolfes]{%
  Friedrich-Alexander-Universität Erlangen-Nürnberg,
  Cauerstr. 11,
  91058 Erlangen,
  Germany}
\email{\{dennis.adelhuette, christian.biefel, martina.kuchlbauer, jan.rolfes\}@fau.de}

\newcommand{\CB}[1]{\todo[author=CB,color=green!50,size=\small]{#1}}
\newcommand{\CBil}[1]{\todo[inline,author=CB,color=green!50,size=\small]{#1}}
\newcommand{\MK}[1]{\todo[author=MK,color=yellow!50,size=\small]{#1}}
\newcommand{\MKil}[1]{\todo[inline,author=MK,color=yellow!50,size=\small]{#1}}
\newcommand{\JR}[1]{\todo[author=JR,color=red!50,size=\small]{#1}}
\newcommand{\JRil}[1]{\todo[inline,author=JR,color=red!50,size=\small]{#1}}
\newcommand{\DA}[1]{\todo[author=DA,color=orange!50,size=\small]{#1}}
\newcommand{\DAil}[1]{\todo[inline,author=DA,color=orange!50,size=\small]{#1}}

\date{\today}

\begin{abstract}
  Pareto efficiency for robust linear programs was introduced by Iancu and Trichakis in \cite{Iancu2014a}. 
We generalize their approach and theoretical results to robust optimization problems in Euclidean spaces with affine uncertainty. 
Additionally, we demonstrate the value of this approach in an exemplary manner in the area of robust semidefinite programming (SDP). In particular, we prove that computing a Pareto robustly optimal solution for a robust SDP is tractable and illustrate the benefit of such solutions at the example of the maximal eigenvalue problem. 
Furthermore, we modify the famous algorithm of Goemans and Williamson \cite{Goemans1995a} in order to compute cuts for the robust \textsc{max-cut} problem that yield an improved approximation guarantee in non-worst-case scenarios. 
\end{abstract}

\keywords{Semidefinite Programming,
Pareto Optimality,
Robust Optimization%
%
}

\makeatletter
\@namedef{subjclassname@2020}{%
	\textup{2020} Mathematics Subject Classification}
\makeatother

\subjclass[2020]{
%
%
90C17, 
90C22 
}

\maketitle

\section{Introduction}
\label{sec:introduction}
Pareto efficiency is a well-established concept in a variety of fields such as economy, engineering and biology, see e.g. \cite{Stewart2008} for a broad overview. 
In \cite{Iancu2014a}, Iancu and Trichakis adapted this concept to robust optimization (RO) for linear programs. 
In particular, they consider the robust linear program
\begin{equation}\label{Eq: RLP}
  \max_{x\in \X}\min_{p\in \Unc}p^\top x,
\end{equation}
where the feasible set $\X$ and the uncertainty set $\Unc$ are assumed to be polytopes. In this setting they characterize and compute so-called Pareto robustly optimal or PRO solutions. These are robustly optimal solutions $x\in \X$ for which there exists no $\bar{x}\in \X$ such that $p^\top \bar{x} \geq p^\top x$ for all $p\in\Unc$ and $\bar{p}^\top \bar{x} > \bar{p}^\top x$ for at least one $\bar{p} \in \Unc$.
The main purpose of this article is to generalize this definition and retrieve a characterization of PRO solutions in a setting that is similar to the one in \cite{Iancu2014a}. Moreover, we show that in the case of robust semidefinite programs, computing PRO solutions is tractable.

Although the work of Iancu and Trichakis on the linear framework is rather new, it has triggered further research such as an analysis for adjustable settings, see e.g. \cite{Ruiter2016} for a rolling horizon approach and \cite{Bertsimas2020} for a Fourier-Motzkin Elimination based approach. 

\subsection*{Structure}
In Section \ref{Sec:sdp-formulation}, we generalize the approach of Iancu and Trichakis to $\X$ being a subset of a finite dimensional Euclidean vector space and an uncertain parameter that affects the objective affinely and is contained in a compact, convex uncertainty set $\Unc$. In particular, we provide a characterization of Pareto robustly optimal (PRO) solutions in this broader setting, which is our main result. 
This result enables us to prove the tractability of computing a PRO solution in the case of robust semidefinite programming. 
In Sections~\ref{sec:rob_eigenvalue} and \ref{Sec: Applications}, we illustrate how to compute the robust maximal eigenvalue of a class of matrices and consider a variant of the SDP that is at the core of the Goemans-Williamson Algorithm \cite{Goemans1995a}. 
The PRO solutions of the latter, are then used as an input for the algorithm and improve the computed cuts for the robust \textsc{max-cut} problem.

\subsection*{Notation}
In the remainder of this article, the feasible set $\X$ and the uncertainty set $\Unc$ are contained in finite dimensional Euclidean vector spaces. 
In the present article, we will mostly choose for both spaces the space of real symmetric $n \times n$-matrices $\Scal^n$ equipped with the Frobenius inner product $\langle\cdot,\cdot\rangle$, i.e., $(\Scal^n,\langle\cdot,\cdot\rangle)$. For a positive semidefinite matrix $X \in \R^{n \times n}$, we write $X\succeq 0$ and we denote the set of symmetric positive semidefinite matrices by $\Scal^n_{\succeq 0}$. 
Given a subset $S$ of an Euclidean vector space $V$ with inner product $\langle \cdot, \cdot \rangle _V$, we denote its dual cone by $S^*=\{y\in V\defsep \Mprod{y,x}_V \geq 0~\forall x\in S\}$ and its relative interior by $\mathrm{relint}(S)$. For a real matrix $A \in \R^{n \times n}$, we denote its trace by $\mathrm{Tr}(A)$. For a positive integer $n\in\N$, we use $[n]:=\{1,...,n\}$ to denote a set of indices and $I_n$ to denote the $n$-dimensional identity matrix. The vector $e_i \in \R^n$, $i \in [n]$, denotes the $i$-th unit vector and $\mathbbm{1}:=\sum_{i=1}^n e_i\in \R^n$ denotes the all-ones vector. We further denote by $E_{ij}\coloneqq \frac{1}{2}(e_ie_j^\top+e_je_i^\top)\in \Scal^n$, $i,j \in [n]$, the standard basis of $\Scal^n$.

\section{Pareto optimal solutions for affine uncertainty}\label{Sec:sdp-formulation}
As a generalization of Program~\eqref{Eq: RLP}, we consider the following robust optimization problem
\begin{equation}\label{eq:BasicProblem}
	\sup_{x \in \X} \min_{p \in \Unc} f(x,p),
\end{equation}
where $\X$ is the feasible set, $\Unc\subseteq V$ is the convex and compact uncertainty set located in a Euclidean vector space. Let further $f(\cdot,p)\colon \X \to \R$ be a function that is well-defined for all $p \in \Unc$. Naturally, we assume that $\Unc$ is not a singleton. The parameter $p \in \Unc$ encodes an affine uncertainty, i.e., $f(x,\cdot)\colon \Unc \to \R$ is affine in $p$ for all $x \in \mathcal{X}$. The involved affinity gives rise to an alternative formulation of \eqref{eq:BasicProblem}, namely
\begin{equation}\label{eq:BasicProblem_Euclidean}
	\sup_{x \in \X} \min_{p \in \Unc} \langle \bar{f}(x), p\rangle_V + g(x),
\end{equation}
where $\bar{f}(x)\in V$ and $g(x) \in \R$ are the unique elements that correspond to the affine functional $f(x,\cdot): p\mapsto f(x,p)$ as given by the Riesz' representation theorem. Hence \eqref{eq:BasicProblem} can be seen as an generalization of \eqref{Eq: RLP} to Euclidean vector spaces. However, over the course of the present article we mainly stick to Formulation~\eqref{eq:BasicProblem}. We note further, that if $\X$ is compact and $f$ is continuous on $\X$, we replace '$\sup$' by '$\max$' in \eqref{eq:BasicProblem}. 
We denote the set of robustly optimal solutions, i.e. the set of optimal solutions of \eqref{eq:BasicProblem}, by $\XRO$.

In robust optimization, one usually focuses on the worst-case scenario, i.e. it suffices to find any robust solution $x\in \XRO$. In contrast to this approach, we aim for a specific $x\in \XRO$ that also performs well under all other scenarios $p\in \Unc$. 
To this end, we use the definition of Pareto robustness from \cite{Bertsimas2020}, which is a generalization of the definition from \cite{Iancu2014a} as mentioned in the introduction:
\begin{definition}
A robustly optimal solution $x\in \XRO$ is called a \emph{Pareto robustly optimal solution} (PRO) of \eqref{eq:BasicProblem} if there exists no $\bar{x}\in \X$ such that
\begin{align}\label{Eq: Pareto1}
	&\forall p \in \Unc\colon\quad f(\bar{x},p) \geq f(x,p),\\
	&\exists \bar{p} \in \Unc\colon\quad f(\bar{x},\bar{p}) > f(x,\bar{p}).\label{Eq: Pareto2}
\end{align}
In this case, we also write $x\in\XPRO$. 
If $x\notin\XPRO$, we say for an $\bar{x}$, which fulfills \eqref{Eq: Pareto1} and \eqref{Eq: Pareto2}, that it Pareto dominates $x$.
\end{definition}

It is natural to ask whether such solutions exist, if they can be characterized and whether they can be determined properly.

We first give an introductory example that fits into the setting of \eqref{eq:BasicProblem}. 
We thereby demonstrate that the choice of a Pareto optimal solution can significantly improve the objective value.
After proving our main result \Cref{Thm: main} on a characterization of PRO solutions, we apply it to the example.
In Section \ref{Sec: Applications}, a more broad discussion of applications will be done. 
\begin{example}\label{Ex: knapsack}
Consider the robust quadratic knapsack problem:
\begin{align*}
		qkp(R,w,d,\Unc)\coloneqq\sup_{x\in \{0,1\}^n}\min_{p\in \Unc}&\ x^\top R(p)x\\
		\mathrm{s.t.} & \ w^\top x \leq d.
\end{align*}
quadratic knapsack problems arise in various applications. 
For illustrative purposes, we consider an example from \cite{Rhys1970a}, where a logistics company wants to construct hubs, that on the one hand maximize the reward function $x^\top Rx$ but on the other hand are restricted by budgetary constraints $w^\top x\leq d$. 
Here, rewards $R_{ij}$ are paid for shipping a good from hub $i$ to hub $j$ and rewards $R_{ii}, R_{jj}$ are paid for additional services at the hubs $i$ and $j$ if there is a shipping.
Uncertainties in the reward matrix $R$ may for example originate from the type of lorry the company uses.

In the following, we demonstrate that there are PRO solutions $x\in \XPRO$ for quadratic knapsack, that Pareto dominate other robust solutions $x\in\XRO\setminus \XPRO$. Moreover, we show that the improvement in the objective can be significant, if $p$ does not attain its worst-case realization.
As an example, let $w=\mathbbm{1}, d=5$ and $R(p)= \mathbbm{1}\mathbbm{1}^\top+E_{ii}(p_1-1)+E_{jj}(p_1-1)+E_{ij}(p_2-1)$ for a fixed pair of indices $i,j\in [n]$. This affine relation is a common form to formulate matrix uncertainties (see e.g. \cite{Ghaoui1998a}). 
It can be generalized by considering arbitrary matrices instead of the standard basis matrices $E_{ij}\in \Scal^n$. 
We consider a convex uncertainty set $\Unc\coloneqq\{p\in \R^2:\ p_1\geq 1, p_1^2\leq p_2, p_2\leq 4\}$ and observe that for this particular $\Unc$ the worst case is attained by $p=(1,1)^\top$ since
\begin{align*}
\min_{p\in \Unc}x^\top R(p)x=\min_{p\in \Unc}(p_1-1)(x_i^2+x_j^2)+(p_2-1)x_ix_j+x^\top \mathbbm{1}\mathbbm{1}^\top x
\end{align*}
and $x\geq 0$.
Hence, in the worst case we have $R(p)=R((1,1)^\top)=\mathbbm{1}\mathbbm{1}^\top$ and consequently every $x\in \{0,1\}^n$ with $\sum_{i \in [n]} x_i =5$ is a robustly optimal solution with objective value $x^\top \mathbbm{1}\mathbbm{1}^\top x=25$.
However, every solution that in addition satisfies $x_i=x_j=1$ Pareto dominates the other robust solutions since the respective objective value is equal to
$$x^\top R(p)x=(p_1-1)(x_i^2+x_j^2)+(p_2-1)x_ix_j + 25=2(p_1-1)+(p_2-1)+25.$$
In our example, the advantage of choosing such an $x\in \XPRO$ compared to a solution $x\in \XRO\setminus \XPRO$ can increase to $30>25$, if $p_1=2$ and $p_2=4$.
\end{example}

The key to characterize and determine PRO solutions is the following theorem which is a generalization of Theorem 1 in \cite{Iancu2014a} and our main result. 

\begin{theorem}\label{Thm: main}
A solution $x^* \in \XRO$ of \eqref{eq:BasicProblem} is PRO if and only if it is an optimal solution to the optimization problem 
\begin{align}\label{Prob:ParetoLinInUnc}
\begin{split}
\sup_{y} & \ f(y,\hat{p})\\
\mathrm{s.t.} & \ \min_{p \in \mathcal{U}} f(y,p) - f(x^*,p) \geq 0, \\
& \ y \in \X 
\end{split}
\end{align}
for an arbitrary $\hat{p} \in \mathrm{relint}(\mathcal{U})$. Every feasible solution $y$ to \eqref{Prob:ParetoLinInUnc} with an objective value greater than $f(x^*,\hat{p})$ Pareto dominates $x^*$. 
Moreover, \rev{if Program~\eqref{eq:BasicProblem} yields an optimal solution then it is PRO.}
\end{theorem} 
\begin{proof}
We begin by pointing out that $\mathrm{relint}(\Unc) \neq \emptyset$ since $\Unc$ is convex. Furthermore, for the inner minimization program, there exists an optimal solution $p^*$ since the objective is affine and $\Unc$ is compact.

If $y$ is feasible for Program~\eqref{Prob:ParetoLinInUnc} with an objective value greater than $f(x^*,\hat{p})$, then the following holds:
\begin{align*}
	&f(y,p) \geq f(x^*,p) \ \forall p \in \Unc,\\
	&f(y,\hat{p}) > f(x^*,\hat{p}).
\end{align*}
In other words, $y$ Pareto dominates $x^*$. 

Next, we show that $x^* \in \XPRO$ if and only if $x^*$ is an optimal solution of Program~\eqref{Prob:ParetoLinInUnc}. However, we have already shown that, if there exists a feasible solution with greater objective value than $x^*$, i.e., if $x^*$ is not optimal for Program~\eqref{Prob:ParetoLinInUnc}, then $x^* \notin \XPRO$. Thus, we only need to show that optimality of $x^*$  for Program~\eqref{Prob:ParetoLinInUnc} implies $x^* \in \XPRO$. 
We assume that $x^*$ is not Pareto robustly optimal. 
Then there exists a solution $y \in \X$ that Pareto dominates $x^*$ and we obtain 
\begin{equation}\label{Eq: helpThm1}
0 < \max_{p \in \mathcal{U}} f(y,p) - f(x^*,p).
\end{equation} 
Since, on the right-hand side of \eqref{Eq: helpThm1}, we optimize an affine function over a convex set $\Unc$, an optimal solution $\bar{p}$ is w.l.o.g. an extreme point of $\Unc$. Additionally, the convexity of $\Unc$ implies that for $\hat{p} \in \mathrm{relint}(\Unc)$, there exist $p\in \mathcal{U}$ and $\varepsilon \in (0,1)$ such that $\hat{p} = \varepsilon \bar{p} + (1- \varepsilon)p$. In particular, we obtain
\begin{gather*}
f(y,\hat{p}) - f(x^*,\hat{p})= \varepsilon (f(y,\bar{p}) - f(x^*,\bar{p}))+(1-\varepsilon)(f(y,p) - f(x^*,p)) > 0,
\end{gather*}
where the inequality follows from the fact that $\bar{p}$ is a maximizer in \eqref{Eq: helpThm1} and that $y$ is a feasible solution of Program~\eqref{Prob:ParetoLinInUnc}. Hence, $x^*$ is not an optimal solution of Program~\eqref{Prob:ParetoLinInUnc} and the claim follows.

For the last claim in \Cref{Thm: main}, assume that $y^*$ is an optimal solution of Program~\eqref{Prob:ParetoLinInUnc}. Assume for contradiction that $y^* \notin \XPRO$. 
Then, there exist $\bar{p} \in \mathcal{U}$ and $z \in \X$ with $f(z,\bar{p}) > f(y^*,\bar{p})$ and $f(z,p) - f(y^*,p) \geq 0$ for all $p \in \mathcal{U}$. 
However, since 
\begin{gather*}
f(z,p) - f(x^*,p) \geq f(z,p) - f(y^*,p) \geq 0 \ \forall p \in \Unc,
\end{gather*}
$z$ is feasible for Program~\eqref{Prob:ParetoLinInUnc}. Furthermore, analogously to before, $f(z,\bar{p}) > f(y^*,\bar{p})$ implies that $f(z,\hat{p}) > f(y^*, \hat{p})$, i.e., the objective value of $z$ is higher than the objective value of $y^*$ -- contradiction to the optimality of $y^*$.
\end{proof}

We observe that since the function $f$ is affine on a convex set $\Unc$, one could reformulate the minimization problem with its dual cone, KKT--conditions or reformulations given in \cite{Ben-Tal2014}. This property would be beneficial to solve Program~\eqref{Prob:ParetoLinInUnc}. In the following, we apply 
\Cref{Thm: main} to the problem given in \Cref{Ex: knapsack}.

\vspace{0.2cm}

\noindent \textbf{Example 1 continued.}
Without loss of generality we set $i=1$ and $j=2$.
We prove that $x^*=(1,1,1,1,1,0,\ldots ,0)^\top$ is a PRO solution to $qkp(R,\mathbbm{1},5,\Unc)$ with 
$R(p)= \mathbbm{1}\mathbbm{1}^\top+E_{11}(p_1-1)+E_{22}(p_1-1)+E_{12}(p_2-1)$ and $\Unc\coloneqq\{p\in \R^2:\ p_1\geq 1, p_1^2\leq p_2, p_2\leq 4\}$. Consider an arbitrary point $\hat{p}\in \text{relint}(\Unc)$. 
Due to Theorem \ref{Thm: main} it suffices to show that $x^*$ is an optimal solution to

\begin{subequations}\label{Prob:Ex1}
	\begin{align}
		\max_{y} & \ y^\top R(\hat{p})y,\\
		\mathrm{s.t.} & \ \min_{p \in \Unc} y^\top R(p)y - (x^*)^\top R(p)x^* \geq 0,\label{Constr: Ex1_min} \\
		& \ y \in \{0,1\}^n,\\
		& \ \mathbbm{1}^\top y \leq 5. \label{Constr: Ex1_limit}
	\end{align}
\end{subequations} 
Here, we can reformulate Constraint~\eqref{Constr: Ex1_min} since
\begin{align*}
	&\min_{p\in\Unc}~(\mathbbm{1}^\top y)^2  -(\mathbbm{1}^\top x^*)^2 +(p_1-1)(y_1^2-(x^*_1)^2)+(p_1-1)(y_2^2-(x^*_2)^2)\\&\quad\quad+(p_2-1)(y_1y_2-x^*_1x^*_2)\\
	 = &\min_{p\in\Unc}~(\mathbbm{1}^\top y)^2 - 25 +(p_1-1)(y_1^2-1)+(p_1-1)(y_2^2-1)+(p_2-1)(y_1y_2-1)\\
	 = &(\mathbbm{1}^\top y)^2 - 25 +(y_1^2-1)+(y_2^2-1)+3(y_1y_2-1),
\end{align*}
where the last equation holds since $p=(2,4)^\top$ is a minimizer for every binary $y$. Moreover, since Constraint~\eqref{Constr: Ex1_limit} implies that $(\mathbbm{1}^\top y)^2 - 25 \leq 0$, we conclude that $y_1=y_2=1$ for every feasible $y \in \{0,1\}^n$. Thus, we reformulate Program~\eqref{Prob:Ex1} to
\begin{subequations}\label{Prob:Ex1_reform}
	\begin{align}
		\max_{y} & \ y^\top R(\hat{p})y\\
		\mathrm{s.t.} & \ y_1= y_2 = 1, \\
		& \ y \in \{0,1\}^n,\\
		& \ \sum_{i=3}^n y_i\leq 3. 
	\end{align}
\end{subequations}
Hence, we have that $y^\top R(\hat{p})y=(x^*)^\top R(\hat{p})x^*$ for all feasible $y$ and conclude that $x^*$ is optimal for \eqref{Prob:Ex1}. 

\vspace{0.25cm}

We observe that the reformulated Program~\eqref{Prob:Ex1_reform} is also a quadratic knapsack problem. Furthermore, the uncertainty set chosen in Example \ref{Ex: knapsack} is an intersection of the second order cone with two halfspaces. We computed a Pareto optimal solution and also checked the Pareto optimality by applying Theorem 1, both by hand. However, an SOCP structure in the uncertainty set as illustrated in the above example may in some cases also allow us to dualize the inner minimization program. Since this dualization approach would result in a convex MINLP even for wider classes of programs under uncertainty, the example suggests that obtaining PRO solutions might be computationally tractable in practice for a variety of problems. However, investigating such properties would be the content of future research.

Another way to determine a PRO solution is given by the following theorem in case one can provide a closed form of $\XRO$:

\begin{theorem}\label{Thm:SubsetXPRO}
Let $\hat{p} \in \mathrm{relint}(\mathcal{U})$. Then $\mathrm{argsup}_{x \in \XRO} f(x,\hat{p})$ is a subset of Pareto robustly optimal solutions of $\eqref{eq:BasicProblem}$.
\end{theorem}

\begin{proof}
Assume that $x^* \in \mathrm{argsup}_{x \in \XRO} f(x,\hat{p})$ but $x^* \notin \XPRO$. 
Then there exists $y\in \XRO$ with $f(x^*,p) \leq f(y,p)$ for all $p \in \mathcal{U}$ and $\bar{p} \in \mathcal{U}$ with $f(x^*,\bar{p}) < f(y,\bar{p})$. Similar to the proof of \Cref{Thm: main}, $\hat{p} = \varepsilon \bar{p} + (1-\varepsilon) p$ for a $p \in \mathcal{U}$ and $\varepsilon \in (0,1)$ holds. Hence,
\begin{gather*}
0 \geq f(y,\hat{p}) - f(x^*,\hat{p}) = \varepsilon (f(y,\bar{p}) - f(x^*,\bar{p})) + (1-\varepsilon) (f(y,p) - f(x^*,p)) > 0,
\end{gather*}
where the first inequality holds since $x^*$ was a maximizer of $f(\cdot,\hat{p})$.
\end{proof}

In contrast to Theorems~\ref{Thm: main} and \ref{Thm:SubsetXPRO}, which aim to determine PRO solutions, the following theorem addresses the question whether there exist non-trivial PRO solutions $x$ for \eqref{eq:BasicProblem}, i.e., $x\in \XPRO$ but $\XPRO \neq \XRO$.

\begin{theorem}\label{Thm: XRO=XPRO?}
Let $\hat{p} \in \mathrm{relint}(\mathcal{U})$ and consider the optimization problem
\begin{align}\label{Prob:ParetoLinInUncII}
\begin{split}
\sup_{x,y} & \ f(y,\hat{p}) - f(x,\hat{p})\\
\mathrm{s.t.} & \ \min_{p \in \mathcal{U}} f(y,p) - f(x,p) \geq 0, \\
& \ y \in \X, \\
& \ x \in \XRO.
\end{split}
\end{align}
Then $\XPRO = \XRO$ if and only if the optimal value of \eqref{Prob:ParetoLinInUncII} equals zero.
\end{theorem}
\begin{proof}
Suppose that there exists a feasible solution $(x^*,y^*)$ of \eqref{Prob:ParetoLinInUncII} with strictly positive objective value. We observe that 
$$\min_{p \in \mathcal{U}} f(y^*,p) - f(x^*,p) \geq 0 \text{ and }f(y^*,\hat{p}) - f(x^*, \hat{p}) > 0$$
implies that $y^*$ Pareto dominates $x^*\in \XRO$ and thus $x^*\in \XRO\setminus \XPRO$.
For the opposite direction, we consider an arbitrary $\bar{x}\in \XRO$ and suppose that the optimal value of \eqref{Prob:ParetoLinInUncII} is zero. This implies that 
\begin{align*}
	\begin{split}
		f(\bar{x},\hat{p})\geq \sup_{y} & \ f(y,\hat{p})\\
		\mathrm{s.t.} & \ \min_{p \in \mathcal{U}} f(y,p) - f(\bar{x},p) \geq 0, \\
		& \ y \in \X.
	\end{split}
\end{align*}
Moreover, equality holds since $y=\bar{x}$ is a feasible and optimal solution and thus we can apply Theorem \ref{Thm: main} to obtain that $\bar{x}\in \XPRO$ and conclude $\XPRO=\XRO$.
\end{proof}

\subsection{A tractable reformulation for SDPs under linear perturbations}\label{Sec: Tractable Pareto Reformulation for SDP with linear Perturbations}
We illustrate the above results by the example of semidefinite programming with uncertainties that solely affect the cost matrix. In addition, we provide a tractability result for this class of optimization problems. We consider a feasible set given by an arbitrary spectrahedron 
\begin{align*}
	\X=\{X\in {\Scal^n_{\succeq 0}}: ~\Mprod{A_j,X} =b_j,~ \forall j\in [k]\},
\end{align*}
and an uncertainty set
\begin{align}\label{eq:matrix-uncertainty}
	\Unc=\left\{P=P_0+\sum_{i=1}^N\mu_iP_i:\ \mu\in [\mu^-,\mu^+]\right\}
\end{align}
with fixed parameters $P_0,\ldots , P_N\in \Scal^n, \mu^-,\mu^+\in \R^N$.
This uncertainty set has been widely used for matrix uncertainty, cf. \cite{Ghaoui1998a}. 
We observe that since the Frobenius inner product $f(X,P)=\Mprod{P,X}$ is bilinear, it encodes linearity in $X$ and in the uncertain parameter $P$. 
Hence, it can be used as an objective function for \eqref{eq:BasicProblem}. 
Thus, we consider the following SDP under cost uncertainty which fits in our setting
\begin{align}\label{unc-prob-1}
\begin{split}
	\sup_{X\in \Scal^n_{\succeq 0}}~\min_{P\in \Unc}~&\Mprod{P,X}\\
	\text{s.t.} ~&\Mprod{A_j,X} =b_j,~~~\forall j\in [k].
\end{split}
\end{align}
It is worth noting that the above problem formulation differs from the more established ones in, e.g. \cite{Ghaoui1998a} or \cite{Ben-Tal1998a} by considering uncertainties in the objective instead of uncertainties in the constraints. Although we do not investigate the exact relation between these two approaches here, we want to point out that the considered problem is a semidefinite version of the setting investigated by \cite{Iancu2014a}. We recall that we aim to compute a Pareto robustly optimal solution for \eqref{unc-prob-1}, i.e., a robustly optimal solution $X\in \XRO$, such that there is no other $\bar{X}\in \X$ that satisfies
\begin{align*}
	& \forall P \in \Unc: \ \Mprod{P,\bar{X}}\geq \Mprod{P,X},\\ 
	& \exists \bar{P} \in \Unc: \ \Mprod{\bar{P},\bar{X}}> \Mprod{\bar{P},X}.
\end{align*}	

The following proposition shows how Theorem \ref{Thm: main} can be used to achieve this. 

\begin{prop}\label{Prop: tractableparetosdp}
A solution $X \in \XRO$ is Pareto robustly optimal for \eqref{unc-prob-1} if and only if the optimal value of
\begin{align}\label{Prob:ParetoSDPI}
\begin{split}
\sup_{Z} & \ \langle\hat{P},Z\rangle\\
\mathrm{s.t.} & \ Z \in \Unc^*, \\
& \ X + Z \in \X 
\end{split}
\end{align}
is $0$. If it is positive with optimal solution $Z$, then $X+Z \in \XPRO$. Moreover, if a PRO solution to \eqref{unc-prob-1} exists,  Program~\eqref{Prob:ParetoSDPI} computes \rev{a PRO solution to \eqref{unc-prob-1}. The corresponding runtime is polynomial in $n$.}
\end{prop}
\begin{proof}
Applying Theorem \ref{Thm: main}, one obtains that $X \in \XRO$ is Pareto robustly optimal if and only if 
\begin{subequations}
	\begin{align}
		\sup_{Y} & \ \langle \hat{P}, Y \rangle,\\
		\mathrm{s.t.} & \ \min_{P \in \mathcal{U}} \langle Y - X, P \rangle \geq 0, \label{Eq:ParetoSDPconstr}\\
		& \ Y \in \X
	\end{align}
\end{subequations}
has an optimal value of $\langle \hat{P}, X \rangle$. Let $Z:= Y-X$. Then, $\langle \hat{P}, Y \rangle \geq \langle \hat{P}, X \rangle$ is equivalent to $\langle \hat{P}, Z \rangle \geq 0$ and the inequality $\min_{P \in \mathcal{U}} \langle Y - X, P \rangle \geq 0$ is equivalent to $Z  \in \Unc^*$, which proves the first part of the claim. In order to prove tractability, we observe 
\begin{align*}
\eqref{Eq:ParetoSDPconstr}  \Leftrightarrow && 0 &\leq \min_{\mu \in [\mu^-,\mu^+]} \langle Y - X, D_0 \rangle + \sum_{i=1}^N \mu_i \langle Y - X, D_i \rangle\\
	 \Leftrightarrow && -\langle Y - X, D_0 \rangle &\leq \min_{\mu \in [\mu^-,\mu^+]} \sum_{i=1}^N\mu_i   \langle Y - X, D_i \rangle\\
	 \Leftrightarrow && -\langle Y - X, D_0 \rangle &\leq \max_{y \in \R^{2n}_{\geq 0}} \left\{ y^\top \begin{pmatrix}
		-\mu^+\\
		\mu^-
	\end{pmatrix}:\ \begin{pmatrix}-I_n & I_n\end{pmatrix}y =\begin{pmatrix}
		\langle Y-X,D_1\rangle \\
		\vdots \\
		\langle Y-X,D_n\rangle
	\end{pmatrix} \right\}
\end{align*}

and consequently, Program~\eqref{Prob:ParetoSDPI} can be written as an SDP which is polynomially solvable in the encoding length of its input: 
\begin{align*}
	\sup_{Y,y} & \ \langle \hat{P}, Y \rangle\\
	\mathrm{s.t.} & \ y^\top \begin{pmatrix}
		-\mu^+\\
		\mu^-
	\end{pmatrix}\geq -\langle Y - X, D_0 \rangle,\\
	& \begin{pmatrix}-I_n & I_n\end{pmatrix}y =\begin{pmatrix}
		\langle Y-X,D_1\rangle \\
		\vdots \\
		\langle Y-X,D_n\rangle
	\end{pmatrix}, \\
	& Y \in \X, y\in \R^{2n}_{\geq 0}.
\end{align*}
We note that this maximization program is computationally tractable since the number of additional variables and constraints is polynomial in the encoding length of the input (namely, $n+1$ additional constraints and $2n$ additional variables). 
\end{proof}
Thus, we have proved that computing a Pareto robustly optimal solution for robust semidefinite programs \eqref{unc-prob-1} with cost uncertainty \eqref{eq:matrix-uncertainty} is tractable. 
In the following section we illustrate its use for a robust eigenvalue problem and the computation of max-cuts on graphs with uncertain weights. 


\section{Application I: The Robust Maximum Eigenvalue Problem}\label{sec:rob_eigenvalue}
In the following paragraphs, we show that computing the maximal eigenvalue of a set of affine combinations of matrices fits into the setting of \eqref{eq:BasicProblem}. The largest eigenvalue problem of a matrix $C$ can be written as (see, e.g., \cite{Overton1993a}):
\begin{equation}\label{Eq: nominalEigenvalue}
\begin{aligned}
\lambda_{\max}&=&\max_{X\in \Scal^n_{\succeq 0}} ~&\Mprod{C,X}&&=&\min_y ~&y\\
&&\text{s.t.}~&\mathrm{Tr}(X)=1~(\Leftrightarrow \Mprod{I_n,X}=1)&&&\text{s.t.}~&yI_n-C\succeq 0.
\end{aligned}
\end{equation}
An optimal matrix $X\in \Scal^n_{\succeq 0}$ for the first optimization problem corresponds to the eigenvector $x$ with respect to the largest eigenvalue $\lambda_{\max}$ of $C$ by $X=xx^\top$. In the remainder of this section, we consider the following robust variant of \eqref{Eq: nominalEigenvalue} with respect to a compact and convex uncertainty set $\Unc$.
\begin{equation}\label{Eq: robustEigenvalueSDP}
	\begin{aligned}
		\lambda_{\max}&=&\max_{X\in \Scal^n_{\succeq 0}}\min_{C\in\Unc} ~&\Mprod{C,X}\\
		&&\text{s.t.}~&\mathrm{Tr}(X)=1.
	\end{aligned}
\end{equation}
Note that for compact and convex uncertainty sets $\Unc$, Sion's minimax theorem \cite{Sion1958a} allows us to interchange the $\max$ and $\min$ operators. Thus, the problem boils down to minimizing the maximal eigenvalue of an affine family of symmetric matrices -- a problem with a wide range of applications, e.g. in stability analysis of dynamic systems or the computation of structured singular values, see \cite{Fan1995a}. In the following example, we provide an instance with non-trivial ($\XPRO\neq \XRO$) Pareto robustly optimal solutions for this eigenvalue problem.
\begin{example}
	Let $C\in \Unc = \left\{\begin{pmatrix}
		1 & 0 \\ 0 & 1
	\end{pmatrix}+\mu\begin{pmatrix}
	1 & -1 \\ -1 & 1
\end{pmatrix}: \mu\in [0,1]\right\}$. Then, the matrix $X'=\frac{1}{2}\begin{pmatrix}
1 & -1 \\ -1 & 1
\end{pmatrix}$ is a robustly optimal solution to \eqref{Eq: robustEigenvalueSDP} since for every $\mu\in [0,1]$ and $X \in \Scal^n_{\succeq 0}$ with $\mathrm{Tr}(X)=1$ we have:
\begin{align*}
	\langle C, X\rangle & = \langle I_2,X\rangle +\mu \left\langle\begin{pmatrix}
		1 & -1 \\ -1 & 1
	\end{pmatrix},X\right\rangle \geq \langle I_2,X\rangle = 1.
\end{align*}
Note that the inequality holds because the matrix $\begin{pmatrix}
1 & -1 \\ -1 & 1
\end{pmatrix}$ is positive semidefinite.
Thus, for every feasible $X$, $\mu=0$ is the worst case realization of uncertainty that can occur. Consequently, every feasible solution $X$, such as $X'$, is also a robustly optimal solution.
However, $X'$ Pareto dominates every other solution $X\in\XRO$, since for every $\mu>0$ and $X\neq X'$, we have
$$\langle C, X\rangle = \langle I_2,X\rangle +\mu \left\langle\begin{pmatrix}
	1 & -1 \\ -1 & 1
\end{pmatrix},X\right\rangle < 1 + \mu\left\langle\begin{pmatrix}
1 & -1 \\ -1 & 1
\end{pmatrix},X'\right\rangle = \langle C, X'\rangle.$$  
We note that one could check $X' \in \XPRO$ by an application of Proposition~\ref{Prop: tractableparetosdp}.
\end{example}

Note that the existence of more than one robustly optimal solution is non-trivial as for uncorrelated uncertainties, i.e. uncorrelated uncertainty sets for the entries of $C$, we often obtain a unique robustly optimal solution. In the above example, the uncertainties in the entries are linked through the matrix $\begin{pmatrix}
	1 & -1 \\ -1 & 1
\end{pmatrix}$ and thus correlated.

\section{Application II: Robust Max-Cut}\label{Sec: Applications}

The \textsc{weighted max-cut} problem is one of the fundamental combinatorial problems from Karp's list of 21 NP-complete problems \cite{Karp1972a}. 
Given an undirected graph $G=(V,E)$ equipped with a weight function $w: E\rightarrow \R$, the task is to find a cut $\delta(V')=\{e\in E: |e\cap V'|=1\}$ defined by $V'\subseteq V$ with maximal weight, i.e.,
\begin{equation*}
	mc(G,w):=\max_{V'\subseteq V} \sum_{e\in \delta(V')}w_e = \max_{x\in \{-1,1\}^V} \frac{1}{4} x^\top L_w x,	
\end{equation*}
where $L_w$ denotes the weighted Laplacian of the graph, i.e. $$L_w=\sum_{\{i,j\}\in E}w_{ij}E_{ij}'\text{ with }E_{ij}'=E_{ii}+E_{jj}-2E_{ij}.$$
In combinatorial optimization under uncertainty, it is common to restrict oneself to uncertainties in the objective in order to keep the structure of the underlying combinatorial problem, see \cite{Kasperski} for a survey. 
In the remainder of this section, we consider uncertain weights, i.e., $w\in\Zc\subseteq \R^E$ for a convex and compact uncertainty set $\Zc$. Similar to \cite{Lasserre2006a}, we define the robust counterpart of the uncertain \textsc{weigthed max-cut} problem that corresponds to $mc(G,w)$ by
\begin{equation}\label{Eq: robustmaxcut}
	mc(G,\Zc)= \max_{x\in \{-1,1\}^V}\min_{w\in \Zc} \frac{1}{4} x^\top L(w) x,
\end{equation}
where $L(w)=\sum_{\{i,j\}\in E}w_{ij}E_{ij}'$ denotes the uncertain Laplacian.
Note that the set $\Unc=\{L(w): w\in \Zc\}$ 
represents a more general uncertainty compared to \eqref{eq:matrix-uncertainty} in the previous section. 
Again, we address the question whether for a given graph $G$, we can improve a robustly optimal solution to \eqref{Eq: robustmaxcut} in terms of Pareto dominance. 
In some instances such as $\gamma$-\emph{stable} graphs introduced by Bilu and Linial \cite{Bilu2010a}, there exist solutions $\hat{x}$ that are not only Pareto optimal but moreover ensures that there is no solution $\bar{x}\in\X$ such that there exists $\bar{p}\in \Unc: f(\bar{x},\bar{p}) > f(\hat{x},\bar{p})$.
Although our techniques would apply for their instances there are more efficient ways to compute these solutions. 
However, in general, graphs are not $\gamma$-stable and hence we first demonstrate the existence of two optimal solutions to an instance of robust weighted \textsc{max-cut} problem of which one Pareto dominates the other with the following example:

\begin{example}\label{Ex:Max-CutToy}
	Consider the complete graph with three nodes equipped with uncertain weights $w_{12}(\mu)=w_{13}(\mu)=4+2\mu$ and $w_{23}(\mu)=3+\mu$ that affinely depend on $\mu$ with $\mu\in [-1,1]$. We observe that 
	$$8+4\mu=w(\delta(v_1))\geq w(\delta(v_2))=w(\delta(v_3))=7+3\mu,$$
	where equality holds if and only if $\mu=-1$. Since this describes the worst case for all these three cuts, we have that every cut is a robustly optimal solution. However, the cut $\delta(v_1)$ Pareto dominates the other cuts, since $w(\delta(v_1))>w(\delta(v_2))=w(\delta(v_3))$ whenever $\mu>-1$. 
\end{example}
	
Additionally to \Cref{Ex:Max-CutToy}, we briefly discuss pure interval uncertainty sets which are commonly used for combinatorial optimization under uncertainty, cf. \cite{Kasperski} and \cite{Buchheim2018}. The following shows that in this case Pareto dominance between robustly optimal solutions is only possible under very specific conditions.

\begin{prop}\label{prop:comb}
Consider Program~\eqref{eq:BasicProblem} with $\X\subseteq \{0,1\}^n$, $f(x,p)=p^\top x$, interval uncertainty $\mathcal{U}:=[\bar{p} - \Delta p, \bar{p}] \subseteq \R^n$, and let $x^* \in \XRO$. Then, $x^* + z$ with $z \in \{-1,0,1\}^n$ Pareto dominates $x^*$ if and only if 
\begin{itemize}
  \item $x^*+z \in \XRO$,
  \item \rev{$\{i\in[n]\defsep z_i=-1\}\subseteq\{i\in[n]\defsep\Delta p_i = 0\}$, and,}
  \item there exists at least one $i \in [n]$ with $z_i = 1$ and $\Delta p_i > 0$.
\end{itemize}
\end{prop}
\begin{proof}
Theorem 1 in \cite{Iancu2014a}, which in this case is 
equivalent to our \Cref{Thm: main}, states that $x^*\in \XRO$ is Pareto dominated by $x^* + z^*$ if and only if, for an arbitrary
$\hat{p} \in \mathrm{relint}(\Unc)$, $z^*$ is feasible to the program
\begin{align}\label{Prob:ParetoLinInUncIII}
  \begin{split}
  \max_{z} & \ \hat{p}^\top z\\
  \mathrm{s.t.} & \ z \in \Unc^*, \\
  & x^*+z \in \X,
  \end{split}
\end{align}
and its objective value is positive. We determine the dual cone:
\begin{align}
z \in \Unc^* & \Leftrightarrow z^\top u \geq 0 \ \forall u \in \Unc, \notag \\
& \Leftrightarrow \min_{u \in [\bar{p} - \Delta p, \bar{p}]} z^\top u \geq 0, \notag\\
& \Leftrightarrow \max_{(y,s)\in\R_{\geq 0}^{2n}: \ y - s = z} (\bar{p} - \Delta p)^\top y - \bar{p}^\top s \geq 0, \label{Eq: help_strong_dual}\\
& \Leftrightarrow \exists y \geq 0: \bar{p}^\top z - \Delta p^\top y \geq 0, \ y \geq z,\notag  
\end{align}
where we apply strong duality to obtain \eqref{Eq: help_strong_dual}.
Since for all $i \in [n]$, there exists $\lambda_i \in (0,1)$, such that $\hat{p}_i = \bar{p}_i - \lambda_i\Delta p_i$, Program~\eqref{Prob:ParetoLinInUncIII} is equivalent to
\begin{align}\label{Prob:ParetoLinInUncIV}
\begin{split}
\max_{y,z} & \ \sum_{i\in [n]} (\bar{p}_i - \lambda_i\Delta p_i)z_i\\
\mathrm{s.t.} & \ \bar{p}^\top z - \Delta p^\top y \geq 0, \\
& x^*+z \in \X, \\
& y \geq z, \\
& y \geq 0.
\end{split}
\end{align}
for $\lambda \in (0,1)^n$. Now, $x^* + z^*$ Pareto dominates $x^*$ if and only if there exists $y^* \in \R^n$ such that $(y^*,z^*)$ is a feasible solution to Program~\eqref{Prob:ParetoLinInUncIV} with positive objective value. Since $\lambda$ is arbitrary, this holds for every $\lambda \in (0,1)^n$. Using this property, we prove the proposition in the following.

We assume that $x^* + z^*$ Pareto dominates $x^*$. Thus, $x^* + z^*, x^* \in \XRO$ and, in particular,
\begin{gather}\label{Eq:help0}
\min_{p \in \Unc} p^\top (x^* + z^*) = \min_{p \in \Unc} p^\top x^*.
\end{gather}
Since $x^*$, and $x^*+z^*$ are nonnegative, the worst-case uncertainty is attained at $\bar{p} - \Delta p$. We obtain $(\bar{p} - \Delta p)^\top (x^* + z^*) = (\bar{p} - \Delta p)^\top x^*$, implying $\bar{p}^\top z^* = \Delta p^\top z^*$. Thus, we can set $y_i = |z^*_i|$, $i \in [n]$, and $z= z^*$ to obtain a feasible solution to \eqref{Prob:ParetoLinInUncIV} with objective value 
\begin{gather}\label{Eq:help}
\sum_{i \in [n]} (1 - \lambda_i) \Delta p_i z^*_i 
\end{gather}
which is strictly positive for every $\lambda \in (0,1)^n$ by Theorem~\ref{Thm: main}. This implies 
\begin{gather}\label{Eq:help2}
\sum_{i \in [n]} (1 - \lambda_i) \Delta p_i z^*_i \geq 0
\end{gather}
for all $\lambda \in [0,1]^n$. 
Thus, Inequality~\eqref{Eq:help2} is also true for $\lambda = \sum_{j \in [n]\setminus \{i\}} e_j$ for all $i \in [n]$.
This implies $\Delta p_i z^*_i \geq 0$ for all $i \in [n]$ and thus, whenever $z^*_i = -1$, $\Delta p_i = 0$. Furthermore, \eqref{Eq:help} can only be positive when there exists an index $i \in [n]$ with $z_i = 1$ and $\Delta p_i > 0$. 

Proving the other direction is rather direct, since $x^* + z^* \in \XRO$ implies Equation~\eqref{Eq:help0} and $(y,z)$ with $y_i = |z^*_i|$, $i \in [n]$, and $z= z^*$ is again a feasible solution to Program~\eqref{Prob:ParetoLinInUncIV}. Since the resp. objective value is strictly positive for $\lambda \in (0,1)^n$, $x^* + z^*$ Pareto dominates $x^*$.
\end{proof}

We observe that $x'\in\XRO$ Pareto dominates $x\in\XRO$ only if there exists at least one index $i\in[n]$ with $x'_i=1$, $x_i=0$ and $\Delta p_i>0$, i.e., 
\rev{there} is a scenario $p \in \Unc$ with $p_i>\bar{p}_i-\Delta p_i$ and 
$p_j=\bar{p}_j-\Delta p_j$ for all $j\neq i$
increasing only the solution $x'$ compared to the worst case.  
Additionally, all indices $i\in[n]$ with $x_i=1$ and $x'_i=0$ cannot be affected by uncertainty. 
This second observation leads to the following corollary.

\begin{coroll}
Consider the setting of \Cref{prop:comb}. 
If $\Delta p > 0$, a solution $x\in \XRO$ \rev{is} Pareto dominated by another solution $x'\in\rev{\X}$ \rev{if} and only if 

\begin{itemize}
\item 
\rev{$\{i\in[n]\defsep x_i=1\}\subsetneq \{i\in[n]\defsep x'_i=1\}$, and,}
\item \rev{$\Delta p_j=\bar{p}_j \fa j\in \{i\in[n]\defsep x'_i=1\}\setminus \{i\in[n]\defsep x_i=1\}$.}
\end{itemize}
\rev{If, in addition to $\Delta p > 0$, $\Delta p_i\neq \bar{p}_i$ for all $i\in[n]$, $\XRO = \XPRO$.}

\end{coroll}

Since \textsc{max-cut} can be phrased as a binary program by using the cut polytope, the statements above hold true for the robust \textsc{max-cut} problem for uncorrelated uncertainties.
Although the nominal \textsc{max-cut} problem is widely considered in the literature, its robust counterpart is to the best of our knowledge not well-investigated. For the nominal case, the famous algorithm of Goemans and Williamson \cite{Goemans1995a} enables us to compute a cut that satisfies an $\alpha$-approximation ratio with $\alpha=0.878...$. Moreover, if Khot's  unique games conjecture \cite{Khot2002a} holds, this is the best approximation ratio we could hope to achieve with a polynomial time algorithm.
In the remainder of this section, we first derive 
robustly optimal cuts with the same approximation ratio and then apply our results from Section \ref{Sec:sdp-formulation} to compute new cuts with improved approximation guarantees if the worst-case uncertainty is not attained. To this end, we consider the SDP relaxation of \eqref{Eq: robustmaxcut}:
\begin{equation}\label{eq:robust_sdp_mc_lw}
	\begin{split}
		sdp(G,\Zc)=\max_{Y \in \Scal^n_{\succeq 0}} \min_{w\in \Zc} \quad&\left\langle \frac{1}{4}L(w),Y\right\rangle\\
		& \mathrm{s.t.} \ \langle E_{ii}, Y\rangle = 1 \quad \forall i\in [n].\\
	\end{split}
\end{equation}
If the inner problem in \eqref{eq:robust_sdp_mc_lw} is a tractable conic program, such as an LP or SDP, it can often be dualized and we can properly compute a robustly optimal solution to \eqref{eq:robust_sdp_mc_lw} by solving the resulting SDP. This solution could then be used to compute a cut via Goemans-Williamson's Algorithm that guarantees the same approximation ratio for the robust \textsc{max-cut}.
\begin{prop}
	Let $w\geq 0$ for every $w\in \Zc$ and $\bar{Y}$ be a robust 
	optimal solution to \eqref{eq:robust_sdp_mc_lw}. Then,
	\begin{align*}
	  \min_{w\in\Zc} \left\langle \frac{L(w)}{4},\bar{Y}\right\rangle 
	    = sdp(G,\Zc)\geq mc(G,\Zc) \geq 0.878\ldots sdp(G,\Zc).
	\end{align*}
\end{prop}

\begin{proof}
The first inequality follows by a simple relaxation argument. 
For the second inequality we strictly follow  the arguments of Goemans and 
Williamson \cite{Goemans1995a}: 

Let $\bar{y_k}$ denote the columns of the Cholesky decomposition of $\bar{Y}$. 
Then, we observe that $x\in \{-1,1\}^V$ defined by 
$x_k=\text{sign}(\bar{y_k}^\top r)$ forms a cut in $G$.
The proof of Goemans and Williamson then relies on the fact that for vectors 
$r\in S^{n-1}$ drawn from the rotationally invariant probability distribution on the unit sphere and their corresponding cuts, we have that
\begin{align*}
	\mathbb{E}\left(1-x_ix_j\right)\geq 0.878\ldots (1-\bar{y_i}^\top \bar{y_j})
	=0.878\ldots sdp(G,\Zc).
\end{align*}
Finally, we conclude 
\begin{align*}
	\mathbb{E}\left(\min_{w\in\Zc} \frac{1}{4} x^\top L(w) x\right)
	& = \mathbb{E}\left(\min_{w\in\Zc}\frac{1}{4}\sum_{\{i,j\}\in E}
	w_{ij}(1-x_ix_j)\right)\\ 
	& = \min_{w\in\Zc}\frac{1}{4}\sum_{\{i,j\}\in E}
	w_{ij}\mathbb{E}\left(1-x_ix_j\right)\\
	& \geq 0.878\ldots \min_{w\in\Zc}\frac{1}{4}\sum_{\{i,j\}\in E}w_{ij}
	(1-\bar{y_i}^\top \bar{y_j})\\
	& = 0.878\ldots sdp(G,\Zc).
\end{align*}		
\end{proof}
It is worth noting that there are already similar approximation results known, see e.g. \cite{Kasperski}. 
We observe that the quality of a cut in a graph with uncertain edge weights may not only rely on its performance in a worst case scenario but also on its performance in every other scenario $w\in\Zc$. Hence, we show that a Pareto optimal solution $Y^*$ to \eqref{Eq: robustmaxcut} outperforms any other robustly optimal solution $\bar{Y}$ of $sdp(G,\Zc)$ in terms of the approximation ratio of their corresponding cuts:
\begin{prop} 
	Let $Y^*$ Pareto dominate $\bar{Y}$ for \eqref{eq:robust_sdp_mc_lw} 
	and let $x^*$ and $\bar{x}$ denote the corresponding cuts derived from $Y^*$ and $\bar{Y}$ 
	respectively via the Goemans-Williamson Algorithm. Denote 
	$$sdp(G,w,Y)=\frac{1}{4}\sum_{\{i,j\}\in E}w_{ij}(1-y_i^\top y_j).$$
	Then, for every $w\in\Zc$ we have
	\begin{align*}
	  mc(G,w)\geq 0.878...sdp(G,w,Y^*) \geq 0.878...sdp(G,w,\bar{Y})
	\end{align*}
	and there exists a $w\in\Zc$, for which the last inequality holds strictly.
\end{prop}
\begin{proof}
	\begin{align*}
		\mathbb{E}\left(\frac{1}{4}\sum_{\{i,j\}\in E}w_{ij}(1-x_ix_j)\right) 
		  & = \frac{1}{4}\sum_{\{i,j\}\in E}w_{ij}\mathbb{E}\left(1-x_ix_j\right)\\
		& \geq 0.878\ldots \frac{1}{4}\sum_{\{i,j\}\in E}w_{ij}(1-(y_i^*)^\top y_j^*)\\
		& \geq 0.878\ldots \frac{1}{4}\sum_{\{i,j\}\in E}w_{ij}(1-\bar{y_i}^\top \bar{y_j}),
	\end{align*}	
where the last inequality and its strict counterpart for at least one realization of the uncertain parameter follows from the Pareto dominance of $Y^*$.
\end{proof}

\section{Conclusion}
\label{sec:conclusion}
In this paper, we generalized the methods introduced in \cite{Iancu2014a} to determine Pareto robustly optimal solutions for linear programs with an uncertain objective to general optimization problems whose objective function is affected affinely by the uncertainty. Moreover, we proved the tractability of these methods in the case of semidefinite programming with matrix box uncertainties and illustrated their use at the examples of the maximal eigenvalue of an affine set of matrices and the classical \textsc{max-cut} problem.

\section*{Data availability}
Data sharing not applicable to this article as no datasets were generated or analysed during the current study.
\section*{Acknowledgments}
\label{sec:acknowledgements}

This research has been supported by funding of the Bavarian State Government within the Energie Campus Nürnberg (EnCN).
%
The authors thank the Deutsche Forschungsgemeinschaft for their
support within project B06 in the
Sonderforschungsbereich/Transregio 154 \enquote{Mathematical
  Modelling, Simulation and Optimization using the Example of Gas
  Networks}. Furthermore, this paper has received funding from the European Union’s Horizon 2020 research and innovation program under the Marie Skłodowska-Curie grant agreement No. 764759. 

\printbibliography
\end{document}